\documentclass[11pt]{article}
\usepackage{amssymb}

\ifx\pdfoutput\undefined
 \usepackage[dvips]{graphicx}
 \else
 \usepackage[pdftex]{graphicx}
 \fi
\usepackage{epstopdf}
\usepackage[dvips,letterpaper,margin=1in]{geometry}
\usepackage{subfigure}

\newtheorem{theorem}{Theorem}[section]
\newtheorem{corollary}[theorem]{Corollary}
\newtheorem{lemma}[theorem]{Lemma}

\newcommand{\sol}{\,\rm{sol}}

\newcommand{\pqr}
{\ensuremath{{\mathfrak A}_p{\mathfrak A}_q{\mathfrak A}_r}}

\newcommand{\fvp}{\ensuremath{f_{\mathfrak V}(p^{\alpha}q^{\beta}r^{\gamma})}} 

\newcommand{\fxp}{\ensuremath{f_{\mathfrak X}(p^{{\alpha}_1}q^{{\beta}_1}r^{{\gamma}_1})}}

\begin{document}

\title{{\itshape Enumeration of Groups in Varieties of $A$-groups: A survey}
\author{{\bf{Geetha Venkataraman}}$^{\rm a}$$^{\ast}$\thanks{$^\ast$Corresponding author. Email: geetha@aud.ac.in
\vspace{6pt}}\\\vspace{6pt}  $^{a}${\em{Ambedkar University Delhi, Lothian Road, Kashmere Gate, Delhi 110006, India.}}\\}}

\maketitle

\begin{abstract}
Let $\cal S$ be a class of groups and let $f_{\cal S}(n)$ be the number of isomorphism classes of groups in $\cal S$ of order $n$. Let $f(n)$ count the number of groups of order $n$ up to isomorphism.  The asymptotic bounds for $f(n)$ behave differently when restricted to abelian groups, $A$-groups and groups in general. We survey some results and some open questions in enumeration of finite groups with a focus on enumerating within varieties of $A$-groups. 

``{\bf Accepted for publication in the Special issue of the Proceedings of Telangana Academy of Sciences}". 

\end{abstract}


{\bf Classcode:} {\bf{20D10, 20C20, 20E10}}

\section{Introduction}

Let $f(n)$ denote the number of isomorphism classes of groups of order $n$. Let $\cal S$ be a class of groups and let $f_{\cal S}(n)$ be the number of isomorphism classes of groups in $\cal S$ of order $n$.  Some interesting classes that have been studied are the class of abelian groups, the class of solvable groups, varieties of groups, $p$-groups and $A$-groups, etc (see \cite{SBPNGV2007}). ($A$-groups are groups whose nilpotent subgroups are abelian.)

The asymptotic bounds for $f(n)$ behave differently when restricted to abelian groups, $A$-groups, and groups in general, primarily due to whether the group itself or its Sylow subgroups are abelian or not. In 1993, L Pyber proved a result which settled a conjecture about $f(n)$. The result \cite{LP1993} was published in The Annals of Mathematics and uses results related to Classification of Simple Finite Groups. Pyber showed that $f(n) \leq n^{(\frac{2}{27} + o(1)){\mu(n)}^2}$ as $\mu(n) \rightarrow \infty$, where $\mu(n)$ represents the highest power to which a prime divides $n$.  

While Pyber's upper bound has the correct leading term, it is certainly not the case for the error term. The key to this puzzle may lie in deeper investigation of $A$-groups. A correct leading term for an upper bound of $f_A(n)$ could lead to the correct error term for an upper bound for $f(n)$. The best we do know about $A$-groups is that $f_{A, \sol}(n) \leq n^{7\mu(n) + 6}$, where $f_{A, \sol}(n)$ is the number of isomorphism classes of solvable $A$-groups of order $n$ (see \cite{GV1997}).  The method to finding the correct upper bound for the enumeration function for $A$-groups may be via enumerating within varieties of $A$-groups. (For varieties of groups see \cite{HN1967}.) 

In this article we survey some of the results in enumeration of finite groups to provide a context for enumeration within a small variety of $A$-groups. The aim of the article is to provide a flavour of the representation theory and counting techniques used in such enumeration problems. Much of the material here is a selection of existing published work. Some unpublished material is also used to illustrate the structure of groups in a particular variety of $A$-groups. 

The article unfolds as follows. The next section surveys main results in enumeration of finite groups with brief commentaries. The third section provides some basic background related to varieties of $A$-groups and specifically the varieties ${\mathfrak U} = {\mathfrak A}_p{\mathfrak A}_q$, and  ${\mathfrak V} = {\mathfrak A}_p{\mathfrak A}_q \vee {\mathfrak A}_q{\mathfrak A}_p$, where $p$ and $q$ are distinct primes.  A discussion on bounds for $f_{\mathfrak V}(n)$ is also presented with a sketch of the main steps of the proof. In the fourth and final section, we present unpublished material related to the structure of a finite group in the variety ${\mathfrak A}_p{\mathfrak A}_q{\mathfrak A}_r$, where $p, q$ and $r$ are distinct primes. The last section will also discuss some open questions related to enumeration in varieties of $A$-groups. Logarithms are taken to the base $2$, unless stated otherwise.

\section{A brief history of group enumerations}

Enumeration questions and questions about classification of groups have been topics of study from even before the time when the notion of an abstract group was defined. From the latter part of the 19th century to the first quarter of the 20th century, several mathematicians worked on classifying and hence enumerating groups of specific types or order. For a good survey of these  results see \cite{HBBEEAB2002}.

In 1895,  Otto H\"{o}lder \cite{OH1895} gave a precise answer to the question about number of groups of order n when n was a product of distinct primes. If $n$ is square-free, that is, a product of distinct primes, then any group $G$ of order $n$ is meta-cyclic and so $G$ will be a semidirect product of a cyclic group by a cyclic group. He was able to use this to prove that if $n$ is square-free then 

$$ f(n) = \sum_{m|n} \prod_{p}   \frac{p^{c(p)}-1}{p-1}$$
 where, in the product, $p$ ranges over prime divisors of $n/m$ and $c(p)$  denotes
 the number of primes $q$ dividing $m$ such that $q \equiv 1 \bmod p$.

A lot of the modern work on group enumerations dates back to a paper \cite{GH1960}, published in 1960, by Graham Higman. In the paper he showed that 
$$
f(p^m) \geq p^{\frac{2}{27}m^3 - O(m^2)} 
$$
where $p$ is prime. In 1965, Charles C Sims \cite{CCS1965} showed that
$$
f(p^m) \leq p^{\frac{2}{27}m^3 + O(m^{\frac{8}{3}})}
.$$ 
An unpublished result of Mike Newman and Craig Seely, referred to and shown in \cite{SBPNGV2007}, brings down the error term to $O(m^\frac{5}{2})$. This led to much speculation on what would be the corresponding estimates for $f(n)$. 

Let $n = p_1^{\alpha_1} \cdots p_k^{\alpha_k}$ be the prime decomposition of $n$. Define $\lambda(n) = \alpha_1 + \cdots + \alpha_k$. Note that when $n = p^m$, $p$ prime, then $\mu(n) = \lambda(n) = m$.

The speculation or conjecture that arose and which was quoted by McIver and Neumann \cite {AMPN1987} was that
$$
f(n) \leq n^{(\frac{2}{27}+ \epsilon) {\lambda(n)}^2}
$$
where $\epsilon \rightarrow 0$ as $\lambda(n) \rightarrow \infty$.

The conjecture stemmed from the feeling that the reason for the number of groups of order $n$ must be because of the large number of choices we have for $p$-groups to chose as Sylow $p$-subgroups rather than the number of ways of putting the groups together.

This sentiment was shown to be right in 1991 by Laszlo Pyber \cite{LP1993}. He showed that the number of groups of order $n$ with specified Sylow subgroups is at most $n^{75\mu + 16}$. This together with the choices available for $p$-groups to be Sylow subgroups, gives the result that $f(n) \leq n^{\frac{2}{27}{\mu(n)}^2 + O({\mu(n)}^{\frac{5}{3}})}$. Again, due to the results for $p$-groups we see that the above upper bound has the right leading term. 
 
 Broadly speaking, the leading term comes from the choices of $p$-groups that are available to be Sylow subgroups and the error term arises from the number of ways in which the Sylow subgroups can be put together to create the required group of order $n$. 
 The error term above is certainly not the best. 
 
 Recall that a finite $A$-group is a group whose Sylow subgroups are abelian. A bound for solvable $A$-groups was first given in 1969 by Gabrielle Dickenson in \cite{GAD1969}. She showed that $f_{A, \sol}(n) \leq n^{c\log n}$ for some constant $c >0$. It was improved by McIver and Neumann in 1987 in \cite{AMPN1987}, where they showed that $f_A(n) \leq n^{\lambda(n) +1}$.  
 Since the number of choices for an abelian group of order $m$ up to isomorphism is at most $m$, Pyber's result on number of groups with specified Sylow subgroups, gives us $f_A(n) \leq n^{75\mu + 17}$. His proof shows that for solvable $A$-groups, we can get $f_{A, \sol}(n) \leq n^{40\mu + 17}$. The best bound know till date is $f_{A, \sol}(n) \leq n^{7\mu(n) + 6}$, shown in \cite{GV1997}. However the bounds in the case of $A$-groups or even solvable $A$-groups is certainly not best possible. So a question still open is that what are the best possible constants $c>0$ and $d>0$ such that $f_A(n) \leq n^{c\mu + d}$?
 
 We can see now that $A$-groups become an important class of groups from the enumeration point of view. Further if a `best possible' value of $c$ is found in the upper bound of $f_A(n)$ then we will be closer to the correct error term in the bound for $f(n)$. Note that, whenever $n$ is a prime power, $f_A(n) \leq n$. It is therefore not possible to hope that there is a constant $c > 0$ such that $n^{c \mu} \leq f_A(n)$ for all $n$. Consequently we need to elaborate further on what is meant by a best possible $c$. We refer to \cite{GV1997} for this. Let ${\cal S}$ be a class of $A$-groups such that there are infinitely many $n$ for which $f_{\cal S}(n) \leq n$. Then $c>0$ will be called best possible if there exists $d > 0$ such that $f_{\cal S}(n) \leq n^{c\mu + d}$ and given any $\epsilon > 0$ there exist infinitely many $n$ with $\mu(n)$ unbounded such that $n^{(c-\epsilon)\mu} < f_{\cal}(n) $.
 
In \cite{AMPN1987} and \cite{LP1993} it was shown that $f_{A,\sol}(n) \geq n^{c\mu}$ for certain specific types of $n$. In both the cases $0 < c < 0.08$.  So the gap between these lower bounds and the best upper bounds we have is huge. The task therefore is to consider certain classes of solvable $A$-groups where the structure of the groups will allow for enumeration techniques leading to `good' upper bounds and by that process also get `good' lower bounds.
 
Let ${\mathfrak A}_n$ denote the variety of abelian groups with exponent $n$ and let ${\mathfrak U} = {\mathfrak A}_p{\mathfrak A}_q$, and  ${\mathfrak V} = {\mathfrak A}_p{\mathfrak A}_q \vee {\mathfrak A}_q{\mathfrak A}_p$, where $p$ and $q$ are distinct primes. In \cite{GV1993}, it was shown that ${\mathfrak U}$ and ${\mathfrak V}$ are both classes of solvable $A$-groups and that there exist positive constants $ c= c_{p, q}, d = d_{p, q}, c'=c'_{p, q}, d' = d'_{p,q}$ such that $f_{\mathfrak U}(n) \leq n^{c\mu + d}$ and $f_{\mathfrak V}(n) \leq n^{c'\mu + d'}$, where $c$ and $c'$ are best possible in the sense discussed above. Further it was shown in \cite{SBPNGV2007} that $c' = \max\{c_{p,q}, c_{q, p}\}$, where $c_{q, p}$ is the leading term for the product variety ${\mathfrak A}_q{\mathfrak A}_p$. The value of $c_{p,q}$ is given by 
$$c_{p,q} = \frac{1}{d} - 2\sqrt {{(\frac{\log p}{\log q})}^2 + \frac{\log p}{d\log q}} + 2 \frac{\log p}{\log q},$$ 
where $d$ is the order of $p$ modulo $q$. In \cite{GV1999}, Sophie Germain primes were considered. These are primes $q$ such that $p=2q+1$ is also prime. For such $p, q$ we see that $d=1$ and so the chance of getting  larger value for $c_{p, q}$ occurs. It is not known if there are infinitely many Sophie Germain primes, but if we assume this and let $q$ tend to infinity, then $\frac{\log p}{\log q} \rightarrow 1$ and so $c_{p, q} \rightarrow 3-2\sqrt2 = 0.17157\ldots$. If $q$ is the largest currently known Sophie Germain prime, namely $q = (2618163402417 \times 2^{1290000}) -1$ and $p =2q+1$, then $c_{p, q}$ already agrees with $3 - 2 \sqrt 2$ to several decimal places. So $c= 0.171$ is the `best' possible lower bound that we have currently for $A$-groups.

\section{Varieties of $A$-groups}
\label{VarAG}

In this section we consider the varieties $\frak U$ and $\frak V$ and give an idea of some of the techniques used to enumerate within these varieties. The original work for this was done in \cite{GV1993} and then presented with modifications in \cite{SBPNGV2007}.

A variety ${\mathfrak D}$ is said to be a variety of $A$-groups if all its nilpotent groups are abelian. It is sufficient to check if the finite nilpotent groups in the variety ${\mathfrak D}$ are abelian. Consider the product variety ${\mathfrak U} = {\mathfrak A}_p{\mathfrak A}_q$, where $p$ and $q$ are distinct primes. Any finite group in this variety is an extension of an elementary abelian $p$-group by an elementary abelian $q$-group and by the Schur-Zassenhaus Theorem it is a semi-direct product.

Let us now consider ${\mathfrak V} = {\mathfrak A}_p{\mathfrak A}_q \vee {\mathfrak A}_q{\mathfrak A}_p$. The finite nilpotent groups in ${\mathfrak V}$ are abelian and any group in ${\mathfrak V}$ is solvable. Indeed ${\mathfrak V}$ is locally finite. Thus ${\mathfrak V}$ is a locally finite solvable variety of $A$-groups. Further it can be shown that  ${\mathfrak V}$ has exponent $pq$, any finite group of order $p^\alpha q^\beta$ and elementary abelian Sylow subgroups lies in ${\mathfrak V}$. The converse is also true. Any finite group $G \in {\mathfrak V}$ has order $p^{\alpha}q^{\beta}$ and its Sylow-subgroups will be elementary abelian.

Next we present a sketch proof of the main steps involved in finding the `best' bounds for ${\mathfrak U}$ and ${\mathfrak V}$.
As mentioned above any finite group in ${\mathfrak U}= {\mathfrak A}_p{\mathfrak A}_q$ is a semi-direct product of its elementary abelian Sylow $p$-subgroup by an elementary abelian Sylow $q$-subgroup.  That is, if $G \in {\mathfrak U}$, then $G = P \rtimes Q$ where $P$ is a Sylow $p$-subgroup of $G$ and $Q$ is a Sylow $q$-subgroup of $G$.

To count the number of isomorphism classes of groups of order $p^{\alpha} q^{\beta}$  in ${\mathfrak U}$, it suffices to count the number of isomorphism classes of groups $P \rtimes_{\theta} Q$ where $P$ is a fixed elementary abelian group of order $p^{\alpha} $, $Q$ is a fixed elementary abelian group of order $q^{\beta}$  and  $\theta$ runs over all homomorphisms from $Q$ to ${\rm Aut}  P$.  Indeed, two possible approaches could be taken at this stage. One would be to regard $\theta(Q)$ as a subgroup of ${\rm Aut}  P$, and count the possibilities. The other would be to regard $P$ as an $\alpha$-dimensional ${\bf F}_pQ$-module and to then count the possibilities for $P$ up to isomorphism. The approach described below, essentially follows the second route. 

We explore the structure of finite groups in ${\mathfrak U}$ further. Define ${\cal X} := \{ G \in {\mathfrak A}_p{\mathfrak A}_q \mid \mbox{$G$ is finite and $Z(G) = 1$}\}$ and let ${\cal Y}$ be the same as ${\cal X}$ with $p$ and $q$ interchanged. Then for any finite group $G$ in ${\mathfrak U}$, there exists a group $G_1 \in {\cal X}$ such that $G = G_1 \times Z(G)$. Further, the isomorphism class of $G_1$ in ${\cal X}$ determines the isomorphism class of $G$ in ${\mathfrak U}$.

Let $G$ be a group in ${\cal X}$. Then $G$ has a normal Sylow $p$-subgroup. Let us denote it by $P$ and let $Q$ be any Sylow $q$-subgroup of $G$. As $Z(G) = 1$, we have $Q$ acting faithfully by conjugation on $P$. Thus $P$ is a ${\bf F}_pQ$-module and it has no non-zero trivial submodules.

Now let $\alpha$ and $\beta$ be natural numbers. Let $Q$ be an elementary abelian $q$-group of order $q^{\beta}$ and let $P$ be a ${\bf F}_pQ$-module of dimension $\alpha$. We shall say that the ${\bf F}_pQ$-module $P$ is of \mbox{type (1)} if $Q$ acts faithfully on $P$ and $P$ has no non-zero trivial ${\bf F}_pQ$-submodule. 

Let $f_{\cal X}(p^{\alpha} q^{\beta})$ denote the number of groups of order $p^{\alpha}q^{\beta}$ in ${\cal X}$ up to isomorphism. Then $f_{\cal X}(p^{\alpha} q^{\beta})$ is the number of orbits under the action of ${\rm Aut} Q$ on the isomorphism classes of ${\bf F}_pQ$-modules of dimension $\alpha$ and \mbox{type (1)}. Note that, by Maschke's Theorem the \mbox{type (1)} modules are completely reducible.

It is shown in \cite{GV1993} that each orbit under the action of ${\rm Aut} Q$ contains another special type of $\alpha$-dimensional ${\bf F}_pQ$-module that are specifically constructed using certain irreducible ${\bf F}_pQ$-modules which have certain chosen subgroups of $Q$ as kernels of the action of $Q$ on these modules. If we call these special representatives as \mbox{type (2)} modules, then  $f_{\cal X}(p^{\alpha} q^{\beta})$ will be bounded above by the number of \mbox{type (2)} modules up to isomorphism.

Using the above  it can be shown that $f_{\mathfrak U}(n) \leq n^{c_{p, q}\mu(n) +1}$. Here 
$$c_{p,q} = \frac{1}{d} - 2\sqrt {{(\frac{\log p}{\log q})}^2 + \frac{\log p}{d\log q}} + 2 \frac{\log p}{\log q},$$ and $d$ is the order of $p$ modulo $q$. 
Further for every $\epsilon > 0$ there exist infinitely many $n$ such that $f_{\mathfrak U}(n) > n^{(c_{p, q}-\epsilon)\mu(n)}$. So $c_{p, q}$ is best possible. 

Let $G$  be a finite group in ${\mathfrak V}= {\mathfrak A}_p{\mathfrak A}_q \vee {\mathfrak A}_q{\mathfrak A}_p$ then $ G \cong X \times Y$ where $X \in {\mathfrak A}_p{\mathfrak A}_q $ and $Y \in {\mathfrak A}_q{\mathfrak A}_p$. 
Using this, we can show that $f_{\mathfrak V}(n) \leq n^{d\mu(n) +2}$ where $d = \max{\{c_{p,q}, c_{q, p}\}}$. 
Further for every $\epsilon > 0$ it can be shown that there exist infinitely many $n$ such that $f_{\mathfrak V}(n) > n^{(d-\epsilon)\mu(n)}$. So $d$ is best possible. 

\section{Structure of groups in ${\mathfrak A}_p{\mathfrak A}_q{\mathfrak A}_r$}

As seen in the earlier discussions we still do not have `best' bounds for $A$-groups or even solvable $A$-groups. While enumerating in small varieties of $A$-groups did gives a class of $A$-groups for which `best' bounds exist, these still do not help us bridge the gap between the upper and lower bounds for $f_{A,\sol}(n)$. The problem with the small variety of $A$-groups considered above is that they did not provide a large enough collection of $A$-groups, up to isomorphism of a given order,  to be able to build a good lower bound. 

To decrease the gap beween $7$ and $3 -2\sqrt2$ we could enumerate in another larger variety of $A$-groups, namely,  ${\mathfrak T} = \vee_{\sigma \in S_3} {\mathfrak A}_{\sigma(p)}{\mathfrak A}_{\sigma(q)}{\mathfrak A}_{\sigma(r)}$ where $p, q, r$ are distinct primes and $S_3$ represents the permutation group on $3$ letters. A first step towards analysing the join variety ${\mathfrak T}$ would be to understand the structure of finite groups in ${\mathfrak W} =\pqr$. The next would be to analyse the structure of groups in ${\mathfrak T}$.
In this section we present some of the unpublished work done in \cite{GV1999} on the structure of groups in  ${\mathfrak W}$. The first lemma shows us that finite groups in ${\mathfrak W}$ are solvable $A$-groups.

\begin{lemma}
\label{pqr_sol}
Let $p$, $q$ and $r$ be distinct primes and let $G$ be a finite
group in ${\mathfrak W} =\pqr$. Then 
\begin{enumerate}
\item[{\rm(}i\rm{)}] $|G| = p^{\alpha}q^{\beta}r^{\gamma}$ for some
$\alpha, \beta , \gamma$ in ${\bf N}$;
\item[{\rm(}ii\rm{)}] $G = PQR$ where $P$ is the {\rm(}unique{\rm)} normal Sylow
$p$-subgroup of $G$; $Q$ and $R$ are some Sylow $q$ and Sylow
$r$-subgroups respectively. Further $P$, $Q$ and $R$ are
elementary abelian groups and $QR \in {\mathfrak A}_q{\mathfrak A}_r$.
\item[{\rm(}iii\rm{)}] $PQ$ is a normal subgroup of $G$.
\item[{\rm(}iv\rm{)}] $G$ is solvable.
\end{enumerate}
\end{lemma}

\noindent{\bf Proof}{ }Product of varieties is associative and so $G \in {\mathfrak A}_p({\mathfrak A}_q{\mathfrak A}_r)$.
Therefore by the Schur-Zassenhaus Theorem, there exists a subgroup $H$ in $G$
such that $G = P \rtimes H $ where $P$ is the normal Sylow $p$-subgroup of $G$. A similar argument shows that $H = Q \rtimes R$ where $Q$ is the normal Sylow $q$-subgroup of $H$ and $R$ a Sylow $r$-subgroup of $H$. Thus $G=PQR$ and $P$, $Q$ and $R$ are elementary abelian $p$, $q$ and $r$ Sylow subgroups of $G$ respectively. Therefore $G$ has the required order. In order to show that $PQ$ is normal in $G$, we note that $R$ normalises both $P$ and $Q$ and hence it normalises $PQ$. For the last part we note that $P$ is abelian and hence solvable. Further $G/P$ is in the metabelian variety ${\mathfrak A}_q{\mathfrak A}_r$ and so is solvable. Therefore $G$ is solvable.  

\begin{theorem}
\label{pqr_sep-cent}
Let $G$ be a finite group in \pqr where $p$, $q$ and $r$ are
distinct primes. Then there exists a $G_0 \in
\pqr$ such that $Z(G_0)$, the centre of $G_0$, is identity and $G =  G_0 \times Z(G)$. Further the choice for $G_0$ is unique up to isomorphism.
\end{theorem} 

\noindent{\bf Proof}{ }By Lemma \ref{pqr_sol} we know that $G=PQR$. We show that there is a subgroup $G_0$, as required, by essentially showing that each Sylow subgroup of $G$ decomposes as a direct product with the corresponding Sylow subgroup of $Z(G)$ being a part.

Now $P$ can be regarded as a ${\bf F}_pQR$-module, where  $p$ does not divide the order of $QR$. So by Maschke's Theorem, $P$ is completely reducible. Let $P_1$ be the Sylow $p$-subgroup of $Z(G)$. Then $P_1$ is a ${\bf F}_pQR$-submodule of $P$. Since $P$ is completely reducible there exists a normal subgroup $P_0$ of $G$ such that $P = P_0 \times P_1$. Further it is obvious that $P_0QR \cap P_1 = \{1\}$. Since $P_1$ is a subgroup of $Z(G)$ every element of $P_1$ commutes with every element of $P_0QR$. Thus $G  = P_0QR \times P_1= {\hat G} \times P_1$ where ${\hat G}= P_0QR$. 

Now $Q$ is a completely reducible ${\bf F}_qR$- module with a submodule $Q_1$ where $Q_1$ is the Sylow $q$-subgroup of $Z(G)$. Thus we can write $Q = Q_0 \times Q_1$ such that $R$ normalises $Q_0$.  By a similar argument as above ${\hat G}  =  P_0Q_0R \times Q_1 = {\bar G} \times Q_1$ where ${\bar G}= P_0Q_0R$. 

For the final step we note that since $R$ is elementary abelian we can write $R = R_0 \times R_1$ where $R_1$ is the Sylow $r$-subgroup of $Z(G)$. Since $P_0Q_0$ is a normal subgroup of ${\bar G}$, we get that ${\bar G} =  P_0Q_0R_0 \times R_1 = G_0 \times R_1$ where $G_0 = P_0Q_0R_0$. Consequently $G = G_0 \times P_1 \times Q_1 \times R_1 = G_0 \times Z(G)$.
Further since $G_0$ is isomorphic to the quotient group, $G/Z(G)$, the choice of $G_0$ is unique up to isomorphism.

\vspace{.1in}
We have the following corollaries to the above theorem. The proof follows obviously from Theorem \ref{pqr_sep-cent}.
\begin{corollary}
Let $p$, $q$ and $r$ be distinct primes. Let 
$$
{\mathfrak X} = \left
\{ X \in \pqr \mid \mbox{$X$ is finite and }Z(X) = 1 \right\}
$$
and let $G$ be a finite group in \pqr.
Then there exists $X \in {\mathfrak X}$ and an abelian group $Z$ such
that $G = X \times Z$. Further if $G = X_1 \times Z_1$ for some
$X_1 \in {\mathfrak X}$ and some abelian group $Z_1$, then $X \cong
X_1$ and $Z \cong Z_1$.
\end{corollary}
\begin{corollary}
\label{pqr_enum}
Let $p$, $q$ and $r$ be distinct primes. Let ${\mathfrak V} = \pqr$.
Then 
$$
\fvp = \sum \fxp
$$ where the sum is over ordered triples of natural numbers
$({\alpha}_1, {\beta}_1, {\gamma}_1)$ such that ${\alpha}_1
\leq \alpha$, ${\beta}_1 \leq \beta$ and ${\gamma}_1 \leq
\gamma$. 
\end{corollary}

We can see from Corollary \ref{pqr_enum} that we now need to investigate the enumeration of groups in ${\mathfrak X}$ of a given order, up to isomorphism. Let $G$ be a group in ${\mathfrak X}$. Then by Lemma \ref{pqr_sol} we know that
$G = PQR$ where $P$ is the normal Sylow $p$-subgroup of $G$. Further $P$ can be regarded as a ${\bf F}_pQR$-module. Since $G$ has a trivial centre, as a ${\bf F}_pQR$-module $P$, has no non-zero trivial submodule. We end this section with a theorem that explains the module structure of $P$ further.

\begin{theorem}
Let $G$ be a group in ${\mathfrak X}$ with $G = PQR$ where $P$, $Q$
and $R$ satisfy the conditions of Lemma \ref{pqr_sol}. Let $H=QR$. Then
\begin{enumerate}
\item[\rm{(}i\rm{)}] $P = P_1 \oplus P_2$, where $P_1$ and $P_2$ are ${\bf
F}_pH$-submodules of $P$ satisfying 
$$
P_1 = \left \{ x \in P \mid h(x) = x \ \mbox{for all}\ h \in Q
 \right \}\: .
$$ 
Further $P_1$ has no non-zero trivial $R$-submodule and $P_2$
does not have any non-zero trivial $Q$-submodule;
\item[\rm{(}ii\rm{)}]If $P'$ is a ${\bf F}_pH$-module such that $P' =
{P_1}' \oplus {P_2}'$ where 
$$
{P_1}'=\left\{x \in P' \mid h(x) = x \ \mbox{for all}\ h \in Q 
\right\}
$$
then $P\cong P'$ as ${\bf F}_pH$-modules if and only if $P_1
\cong {P_1}'$ and ${P_2} \cong {P_2}'$ as ${\bf F}_pH$-modules. 
\end{enumerate}
\end{theorem}

\noindent{\bf Proof}{ } Let $P_1 =C_G(Q) \cap P$. Then $P_1 = O_p(Z(PQ))$. Since $PQ$ is a normal subgroup of $G$, we get that $P_1$ is a normal subgroup of $G$. Thus $P_1$ is a ${\bf F}_pH$-submodule of $P$ and it is obvious that as a ${\bf F}_pH$-module, $P_1 = \left \{ x \in P \mid h(x) = x \ \mbox{for all}\ h \in Q \right \}$.
 
Further, since $p$ does not divide $|H|$, by Maschke's Theorem we have that $P$ is completely reducible. Thus there exists a ${\bf F}_pH$-submodule $P_2$ of $P$ such that $P = P_1 \oplus P_2$. Since we know that $P$ has no non-zero trivial ${\bf
F}_pH$-submodule. Therefore $P_1$ cannot have a non-zero trivial $R$-submodule and $P_2$ cannot have any non-zero trivial $Q$-submodule.

Let $P \cong P'$ as ${\bf F}_pH$-modules via the mapping $\phi$. Then $\phi(P_1) = {P_1}'$. From this we get $P_2 \cong  \phi(P_2) \cong P'/{P_1}' \cong {P_2}'$ as required. The converse is obvious.

\vspace{.05in}
From Corollary \ref{pqr_enum} and the above Theorem, to find a `good' bound for groups in $\pqr$ we need to count in ${\mathfrak X}$. This in turn depends on counting ${\bf F}_pH$-modules $P_1$ and $P_2$ up to isomorphism. Since these are completely reducible, and  $H \in {\mathfrak A}_q{\mathfrak A}_r$, we need to investigate the irreducible ${\bf F}_pH$-modules. A start towards this process was made in \cite{GV2013}. However, much more needs to be done before we are able to get the bounds of the `best' kind.
\vspace{.05in}

We end with a few open questions which were posed in Chapter 22 of \cite{SBPNGV2007}, which are relevant to the discussions in this article. We reproduce them here.
\vspace{.05in}

{\it Question 22.21} Is it the case that $f_A(n) / f_{A,\sol} (n) \rightarrow 1$ as $\lambda(n) \rightarrow \infty$? How big is $f_A(n)  - f_{A,\sol} (n)$  compared with $f_A(n)$? 
\vspace{.05in}

{\it Question 22.22} Define $\alpha = \limsup_{n\rightarrow \infty} \frac{\log f_A(n)}{\mu(n)\log n}$.  What is $\alpha$? Could it perhaps be $3 - 2\sqrt2$?

{\it Question 22.23} For which varieties ${\mathfrak V}$ of $A$-groups is it true that the leading term of the enumeration function $f_{\mathfrak V}(n)$  is equal to the leading term of $f_{\mathfrak U}(n)$  for some minimal non-abelian subvariety ${\mathfrak U}$  of ${\mathfrak V}$?

\bibliographystyle{plain}

\end{document}